\documentclass[A4paper,11pt]{article}
\usepackage{amsfonts}
\usepackage{amsmath}
\usepackage{amssymb}
\usepackage{graphicx}
\usepackage{rlepsf}
\usepackage{amscd}
\begin{document}
\title{New Categorif\mbox{}ications of the Chromatic and the Dichromatic Polynomials for Graphs}
\author{Marko Sto\v si\'c 
\thanks{The author is supported by {\it Funda\c c\~ao para a Ci\^encia e a Tecnologia}/(FCT), grant no. SFRH/BD/6783/2001}\\
Departamento de Matem\'atica  and \\
CEMAT -  Centro de Matem\'atica e Aplica\c c\~oes\\
Instituto Superior T\'ecnico\\
Av. Rovisco Pais 1\\
1049-001 Lisbon\\ 
Portugal\\
e-mail: mstosic@math.ist.utl.pt
}

\date{}

\newtheorem{theorem}{Theorem}
\newtheorem{acknowledgment}[theorem]{Acknowledgment}
\newtheorem{algorithm}[theorem]{Algorithm}
\newtheorem{axiom}[theorem]{axiom}
\newtheorem{case}[theorem]{Case}
\newtheorem{claim}[theorem]{Claim}
\newtheorem{conclusion}[theorem]{Conclusion}
\newtheorem{condition}[theorem]{Condition}
\newtheorem{conjecture}[theorem]{Conjecture}
\newtheorem{corollary}[theorem]{Corollary}
\newtheorem{criterion}[theorem]{Criterion}
\newtheorem{definition}{Def\mbox{}inition}
\newtheorem{example}{Example}
\newtheorem{exercise}[theorem]{Exercise}
\newtheorem{lemma}{\indent Lemma}
\newtheorem{notation}[theorem]{Notation}
\newtheorem{problem}[theorem]{Problem}
\newtheorem{proposition}{Proposition}
\newtheorem{remark}[theorem]{Remark}
\newtheorem{solution}[theorem]{Solution}
\newtheorem{summary}[theorem]{Summary}
\newcommand{\ud}{\mathrm{d}}

\def\gcd{\mathop{\rm gcd}}
\def\Ker{\mathop{\rm Ker}}
\def\max{\mathop{\rm max}}
\def\map{\mathop{\rm map}}
\def\lcm{\mathop{\rm lcm}}
\def\kraj{\hfill\rule{6pt}{6pt}}
\def\diag{\mathop{\rm diag}}
\def\span{\mathop{\rm span}}
\def\deg{\mathop{\rm deg}}
\def\rank{\mathop{\rm rank}}
\def\sgn{\mathop{\rm sgn}}
\def\kvn{\{n\}_q}
\def\F{\mathbb{F}}
\def\R{\mathbb{R}}
\def\C{\mathcal{C}}
\def\P{\mathcal{P}}
\def\A{\mathcal{A}}
\def\D{\mathcal{D}}
\def\H{\mathcal{H}}
\def\N{\mathbb{N}}
\def\K{\mathbb{K}}
\def\Z{\mathbb{Z}}
\def\Q{\mathbb{Q}}
\def\X{\qbezier(0.00,0.00)(0.50,1.00)(1.00,2.00)
\qbezier(1.00,0.00)(0.80,0.40)(0.60,0.80)
\qbezier(0.00,2.00)(0.20,1.60)(0.40,1.20)
}

\def\Y{\qbezier(1.00,0.00)(0.50,1.00)(0.00,2.00)
\qbezier(0.00,0.00)(0.20,0.40)(0.40,0.80)
\qbezier(1.00,2.00)(0.80,1.60)(0.60,1.20)
}

\def\O{\qbezier(0.00,0.00)(0.20,1.00)(0.00,2.00)
\qbezier(1.00,0.00)(0.80,1.00)(1.00,2.00)
}
\def\OF{\qbezier(0.00,0.00)(0.00,1.00)(0.00,2.00)
\qbezier(1.00,0.00)(1.00,1.00)(1.00,2.00)
}

\def\l{
\qbezier(0.00,0.00)(0.50,1.20)(1.00,0.00)
\qbezier(0.00,2.00)(0.50,0.80)(1.00,2.00)
}

\def\LPP{
\qbezier(0.00,0.00)(0.35,0.70)(0.70,1.40)
\qbezier(0.70,0.60)(0.65,0.70)(0.60,0.80)
\qbezier(0.00,2.00)(0.20,1.60)(0.40,1.20)
\qbezier(0.7,1.4)(1.3,2.6)(1.3,1)
\qbezier(0.7,0.6)(1.3,-0.6)(1.3,1)}

\def\LPM{\qbezier(0.00,0.00)(0.2,0.40)(0.40,0.80)
\qbezier(0.60,1.20)(0.65,1.30)(0.70,1.40)
\qbezier(0.00,2.00)(0.35,1.30)(0.70,0.60)
\qbezier(0.7,1.4)(1.3,2.6)(1.3,1)
\qbezier(0.7,0.6)(1.3,-0.6)(1.3,1)}

\def\XS{\qbezier(0.00,0.00)(0.50,1.00)(1.00,2.00)
\qbezier(1.00,0.00)(0.80,0.40)(0.60,0.80)
\qbezier(0.00,2.00)(0.20,1.60)(0.40,1.20)

\qbezier(0.0,2.00)(0.0,1.85)(0,1.70)
\qbezier(1.00,2.00)(1,1.85)(1.0,1.70)
\qbezier(0.00,2.00)(0.10,1.90)(0.20,1.80)
\qbezier(1.00,2.00)(0.9,1.9)(0.8,1.80)}

\def\YS{\qbezier(1.00,0.00)(0.50,1.00)(0.00,2.00)
\qbezier(0.00,0.00)(0.20,0.40)(0.40,0.80)
\qbezier(1.00,2.00)(0.80,1.60)(0.60,1.2)
\qbezier(0.0,2.00)(0.0,1.85)(0,1.70)
\qbezier(1.00,2.00)(1,1.85)(1.0,1.70)
\qbezier(0.00,2.00)(0.10,1.90)(0.20,1.80)
\qbezier(1.00,2.00)(0.9,1.9)(0.8,1.80)}

\def\GR{\qbezier(0.00,0.00)(0.25,0.25)(0.50,0.50)
\qbezier(1.00,0.00)(0.75,0.25)(0.50,0.50)
\qbezier(0.00,2.00)(0.25,1.75)(0.50,1.5)
\qbezier(1.0,2.00)(0.75,1.75)(0.5,1.50)

\qbezier(0.00,2.00)(0,1.9)(0,1.80)
\qbezier(0.00,2.00)(0.1,2)(0.20,2)
\qbezier(1.00,2.00)(1,1.9)(1,1.80)
\qbezier(1.00,2.00)(0.9,2)(0.8,2)
\linethickness{2.5pt}
\qbezier(0.5,0.5)(0.5,1)(0.50,1.5)}

\def\GRF{\qbezier(0.00,0.00)(0.25,0.25)(0.50,0.50)
\qbezier(1.00,0.00)(0.75,0.25)(0.50,0.50)
\qbezier(0.00,2.00)(0.25,1.75)(0.50,1.5)
\qbezier(1.0,2.00)(0.75,1.75)(0.5,1.50)

\linethickness{2.5pt}
\qbezier(0.5,0.5)(0.5,1)(0.50,1.5)}

\def\OS{\qbezier(0.00,0.00)(0.60,1.00)(0.00,2.00)
\qbezier(1.00,0.00)(0.40,1.00)(1.00,2.00)
\qbezier(0.0,2.00)(0.0,1.85)(0,1.70)
\qbezier(1.00,2.00)(1,1.85)(1.0,1.70)
\qbezier(0.00,2.00)(0.10,1.90)(0.20,1.80)
\qbezier(1.00,2.00)(0.9,1.9)(0.8,1.80)}

\maketitle


\begin{abstract}
In this paper, for each graph $G$, we def\mbox{}ine a chain complex of graded modules over the ring of polynomials, whose graded Euler characteristic is equal to the chromatic polynomial of $G$. Furthermore, we def\mbox{}ine a chain complex of doubly-graded modules, whose (doubly) graded Euler characteristic is equal to the dichromatic polynomial of $G$. Both constructions use Koszul complexes, and are similar to the new Khovanov-Rozansky categorif\mbox{}ications of HOMFLYPT polynomial. We also give simplif\mbox{}ied def\mbox{}inition of this triply-graded link homology theory.
\end{abstract}

\section{Introduction}

In \cite{kov} Khovanov introduced the concept of categorif\mbox{}ication of the Jones polynomial for links. For each link $L$ in $S^3$ he def\mbox{}ined a graded chain complex, with grading 
preserving dif\mbox{}ferentials, whose graded Euler characteristic is equal to
the Jones polynomial of the link $L$. This is done by starting from the state sum 
expression for the Jones polynomial (which is written as an alternating sum),
then constructing for each term a module whose graded dimension is 
equal to the value of that term, and f\mbox{}inally, constructing the
dif\mbox{}ferentials 
as appropriate grading preserving maps, so that 
the complex obtained is a link invariant. \\
\indent 
Using similar techniques by starting from the state-sum expression of the chromatic polynomial for graphs, in \cite{gr} was def\mbox{}ined a chain complex of graded modules whose Euler characteristic is chromatic polynomial of a graph. In \cite{moj} was def\mbox{}ined the inf\mbox{}inite series of chain complexes of graded modules (one for each $n\in\N$) whose Euler characteristics are the specializations of the two-variable dichromatic polynomial (and consequently the Tutte polynomial) of graph $G$. The specializations appear since we want to categorify the two-variable polynomial and the ``standard"  techniques of categorifying link (and graph) polynomials (see e.g. \cite{kovroz}, \cite{kov3}, \cite{kol}, \cite{bn1}, \cite{kov}) work only for one-variable polynomials.\\ 
\indent 

In this paper we def\mbox{}ine a chain complex of doubly-graded modules whose doubly-graded Euler characteristic is equal to the whole two-variable dicromatic polynomial.  The idea is partially inspired by the new version of categorif\mbox{}ication of HOMFLYPT polynomial by Khovanov and Rozansky, see \cite{KR2}. They def\mbox{}ined a chain 
 complex of doubly-graded modules whose doubly-graded Euler characteristic is equal to the whole two-variable HOMFLYPT polynomial. In Section \ref{coment} we describe the simplif\mbox{}ied version of their construction (this construction is also implicit in \cite{kovhoh}).\\
\indent Also, we give new categorif\mbox{}ication of the chromatic polynomial for graphs. We do this here in a dif\mbox{}ferent way than in \cite{gr}. We will def\mbox{}ine chain groups (the direct sums of modules corresponding to the vertices of the cube of resolutions) as the cohomologies of certain chain complexes.

\section{Triply graded link homology}\label{coment}

\subsection{Introduction}

\indent In this section we will introduce the parametrization of the HOMFLYPT polynomial that we will categorify. It is very similar to the one in \cite{KR2}. Throughout the chapter we will consider only braid diagrams $D$ of a link $L$, i.e. regular diagrams which are the closures of (upward) oriented braids.\\
\indent  As is well known, every link can be represented by a braid diagram. Also,  the closures of two braid diagrams $D_1$ and $D_2$ are isotopic as oriented links if and only if $D_1$ and $D_2$ are related by a sequence of Markov moves, which are the following (see \cite{mar}):\\
\indent $(i)$ \quad conjugation: $DD'\longleftrightarrow D'D$\\
\indent $(ii)$ \quad transformations in the braid group:
\begin{eqnarray*}
D&\longleftrightarrow& D\sigma_i\sigma_i^{-1}\\
D&\longleftrightarrow& D\sigma_i^{-1}\sigma_i\\
D{\sigma_j}\sigma_i&\longleftrightarrow&D\sigma_i\sigma_j,\quad |i-j|>1\\
D\sigma_i\sigma_{i+1}\sigma_i&\longleftrightarrow&D\sigma_{i+1}\sigma_i\sigma_{i+1}\end{eqnarray*}
\indent $(iii)$\quad transformations $D\longleftrightarrow D\sigma_{n}^{\pm 1}$, for a braid $D$ with $n$ strands.\\
    
\indent In order to def\mbox{}ine the HOMFLYPT polynomial for a link $L$, from its braid diagram representation $D$, we will introduce a function $F$  on braid diagrams with values in the ring of rational functions in $q$ and $t$ def\mbox{}ined uniquely by the following axioms:\\
\indent $\bf{\ast}$  $F(D_1)=F(D_2)$, if $D_1$, $D_2$ are related by Markov move $(i)$\\
\indent $\bf{\ast}$  $F(D_1)=F(D_2)$, if $D_1$, $D_2$ are related by Markov moves $(ii)$\\
\indent $\bf{\ast}$  $F(D\sigma_n)=F(D)$, if a braid $D$ has $n$ strands\\
\indent $\bf{\ast}$  $F(D\sigma_n^{-1})=-t^{-1}q^{-1}F(D)$, if a braid $D$ has $n$ strands\\
\indent $\bf{\ast}$  Skein relation: for every braid diagram $D$ with $n$ strands and $0<i<n$  
$$ q^{-1}F(D\sigma_i)-qF(D\sigma_i^{-1})=(q^{-1}-q)F(D)$$ 
\indent $\bf{\ast}$ If $U$ is the one-strand diagram of the unknot then
$F(U)=1$.\\

\indent In order to obtain a link invariant we need to normalize the function $F$. Let $\alpha=-t^{-1}q^{-1}$ and let 
\begin{equation}
G(D)={\sqrt{\alpha}}^{n_{+}(D)-n_{-}(D)-s(D)+1}F(D),
\label{xxx}
\end{equation}
where $n_+(D)$, $n_-(D)$ and $s(D)$ are the number of positive crossings, negative crossings and the number of strands of $D$, respectively. We denote the expression $n_{+}(D)-n_{-}(D)-s(D)+1$ by $\omega(D)$. Obviously $G(D)$ is invariant under all Markov moves of braids and it satisf\mbox{}ies the HOMFLYPT skein relation
$$ {(q\sqrt{\alpha})^{-1}G(D\sigma_i)-q\sqrt{\alpha}G(D{\sigma_i}^{-1})}=(q^{-1}-q)G(D).$$
Hence, $G(D)$ is equal to the HOMFLYPT polynomial of the link $L$, normalized such that $G(U)=1$.
In Section \ref{kr2a0}, we will def\mbox{}ine a triply graded chain complex $\C(D)$ whose Euler characteristic is equal to $F(D)$.\\
\indent First of all note  our slightly dif\mbox{}ferent convention compared to \cite{KR2} on the value of unknot. This has the advantage that we can  obtain  the  Alexander polynomial directly by specializing $t$ and $q$ ($t=-q$), and the whole sequence of the $n$-specializations of the (reduced) HOMFLYPT polynomial  (see \cite{kovroz}, \cite{moy}). Specif\mbox{}ically, by taking  $t=-q^{1-2n}$ we obtain polynomials $G_n(D)$ that satisfy the skein relation
$$q^{-n}G_n(D\sigma_i)-q^n G_n(D{\sigma_i}^{-1})=(q^{-1}-q)G_n(D), $$
and whose value on the unknot is $G_n(U)=1$. Hence by suitably collapsing the tri-grading to a bi-grading we get a new categorif\mbox{}ication of the $n$-specializations of the HOMFLYPT polynomial. 

\subsection{Graphs with wide edges}\label{grafwide}  

\indent In order to def\mbox{}ine the function $F(D)$ and hence the HOMFLYPT  polynomial $G(D)$, we introduce the trivalent graphs with wide edges as the resolutions of the crossings. Apart from the crossings $\sigma_i$ and $\sigma_i^{-1}$, we introduce wide edges $\bar{E}_i$ placed between the $i$-th and $(i+1)$-th strand of the braid, like in the following picture:

{{
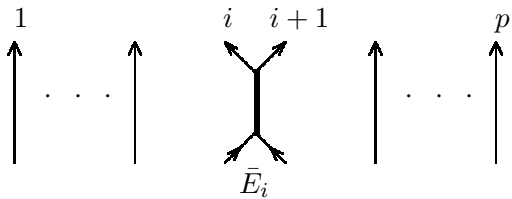
\begin{center}
\setlength{\unitlength}{4mm}
\begin{picture}(30,6) 
\linethickness{0.6pt}
\qbezier(5.00,1.00)(5.00,3.00)(5.00,5.00)
\qbezier(4.80,4.60)(4.90,4.80)(5.00,5.00)
\qbezier(5.20,4.60)(5.10,4.80)(5.00,5.00)

\qbezier(9.00,1.00)(9.00,3.00)(9.00,5.00)
\qbezier(8.80,4.60)(8.90,4.80)(9.00,5.00)
\qbezier(9.20,4.60)(9.10,4.80)(9.00,5.00)

\linethickness{2.0pt}
\qbezier(13.00,2.00)(13.00,3.00)(13.00,4.00)
\linethickness{0.6pt}
\qbezier(12.00,1.00)(12.50,1.50)(13.00,2.00)
\qbezier(14.00,1.00)(13.50,1.50)(13.00,2.00)
\qbezier(12.00,5.00)(12.50,4.50)(13.00,4.00)
\qbezier(14.00,5.00)(13.50,4.50)(13.00,4.00)

\qbezier(14.00,5.00)(13.80,4.90)(13.60,4.80)
\qbezier(14.00,5.00)(13.90,4.80)(13.80,4.60)

\qbezier(12.00,5.00)(12.20,4.90)(12.40,4.80)
\qbezier(12.00,5.00)(12.10,4.80)(12.20,4.60)

\qbezier(12.50,1.50)(12.30,1.40)(12.10,1.30)
\qbezier(12.50,1.50)(12.40,1.30)(12.30,1.10)

\qbezier(13.50,1.50)(13.70,1.40)(13.90,1.30)
\qbezier(13.50,1.50)(13.60,1.30)(13.70,1.10)

\qbezier(17.00,1.00)(17.00,3.00)(17.00,5.00)
\qbezier(16.80,4.60)(16.90,4.80)(17.00,5.00)
\qbezier(17.20,4.60)(17.10,4.80)(17.00,5.00)

\qbezier(21.00,1.00)(21.00,3.00)(21.00,5.00)
\qbezier(20.80,4.60)(20.90,4.80)(21.00,5.00)
\qbezier(21.20,4.60)(21.10,4.80)(21.00,5.00)

\put(6.00,3.00){$\cdot$}
\put(7.00,3.00){$\cdot$}
\put(8.00,3.00){$\cdot$}
\put(18.00,3.00){$\cdot$}
\put(19.00,3.00){$\cdot$}
\put(20.00,3.00){$\cdot$}

\put(5.00,5.50){1}
\put(11.95,5.50){$i$}
\put(13.50,5.50){$i+1$}
\put(21.00,5.50){$p$}
\put(12.50,0.00){$\bar{E}_i$} 
\end{picture}
\end{center}
}}

 Then we can def\mbox{}ine the function $F(D)$ by resolving the crossings by using the following two relations:
\begin{eqnarray}
F(D\sigma_i)&=& F(D\bar{E}_i)-q^2F(D),\label{1}\\
F(D{\sigma_i}^{-1})&=& q^{-2}F(D\bar{E}_i)-q^{-2}F(D).\label{2}
\end{eqnarray}  
Here by $F(D)$ we mean the value of the function $F$ on the diagram that is the closure of the braid diagram $D$, and we have extended the domain of the function $F$ to include trivalent graphs. Then $F$ (restricted to braid diagrams) will satisfy the axioms from the previous subsection, if and only if the values of $F$ on  completely resolved trivalent graphs satisfy
\begin{eqnarray}
&F(U)=1 \label{rel0}\\
&F(D\cup U) = \frac{1+t^{-1}q}{1-q^2} F(D), \quad{\textrm{if }} D {\textrm{ is not an empty diagram}} \label{rel1}\\
&F(D\bar{E}_n)=\frac{1+t^{-1}q^3}{1-q^2} F(D), \quad {\textrm{where }} D {\textrm{  is a diagram with  }} n {\textrm{ strands}}\\
&F(D\bar{E}_i^2)= (1+q^2)F(D\bar{E}_i)\\
&F(D\bar{E}_i\bar{E}_{i+1}\bar{E}_i)+q^2F(D\bar{E}_{i+1})=
F(D\bar{E}_{i+1}\bar{E}_{i}\bar{E}_{i+1})+q^2F(D\bar{E}_{i}).\label{rel4}
\end{eqnarray}

We will use (\ref{1}) and (\ref{2}) in resolving the crossings in the def\mbox{}inition of the triply-graded chain complex that categorif\mbox{}ies  the HOMFLYPT polynomial. 

\subsection{Categorif\mbox{}ication of the two-variable HOMFLYPT polynomial}\label{kr2a0}

In this subsection we will give an alternative, simpler and (essentially) equivalent construction to the one in \cite{KR2} (similar simplif\mbox{}ication is also implicit in \cite{kovhoh}).  \\
\indent Essentially, we will set the variable $a$ from \cite{KR2} to be 0, but we will keep the double grading of the ring of polynomials $R'=\Q[x_1,\ldots,x_{2n}]$. Specif\mbox{}ically, to every arc (line between two crossings) we will assign a dif\mbox{}ferent variable $x_i$, $i=1,\ldots,2n$, where $n$ is the number of crossings of a given planar projection $D$ of a knot $K$. We def\mbox{}ine the bidegree of every $x_i$ to be (0,2) and we put the f\mbox{}ield of coef\mbox{}f\mbox{}icients $\Q$ in bidegree (0,0). Also, by $\{\cdot,\cdot\}$ we denote a shift in bigrading.
\begin{remark}
Note that this corresponds to the case $n=-1$ in \cite{kovroz}, but with the introduction of a new grading direction.   
\end{remark}

Let $L$ be a link and let $D$ be its braid diagram presentation. Let $I$ be the ideal of $R'$ generated by the monomials $x_1+x_2-x_3-x_4$ for every crossing of $D$, and let $R=R'/I$ (see the picture below for the notation).\\
\indent   To each crossing we will assign 0- and 1-resolutions according to the following picture:
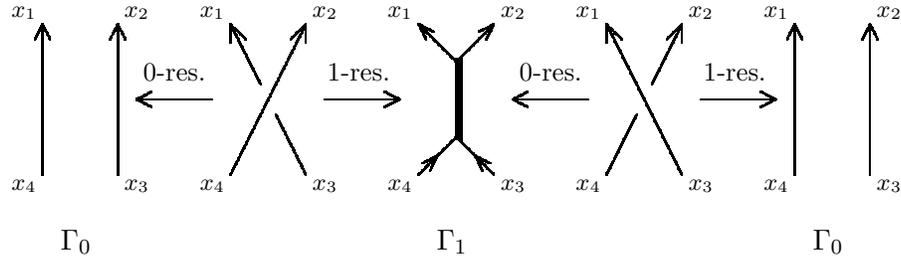
\begin{center}
\setlength{\unitlength}{5mm}
\begin{picture}(25,5) 
\linethickness{0.6pt}
\put(0,-0.5){
\qbezier(4.00,1.00)(4.00,3.00)(4.00,5.00)
\qbezier(3.80,4.60)(3.90,4.80)(4.00,5.00)
\qbezier(4.20,4.60)(4.10,4.80)(4.00,5.00) \qbezier(2.00,1.00)(2.00,3.00)(2.00,5.00)
\qbezier(1.80,4.60)(1.90,4.80)(2.00,5.00)
\qbezier(2.20,4.60)(2.10,4.80)(2.00,5.00) 
\qbezier(9.5,3)(10.5,3)(11.5,3)
\qbezier(11.5,3)(11.3,3.1)(11.1,3.2)
\qbezier(11.5,3)(11.3,2.9)(11.1,2.8)

\qbezier(4.5,3)(5.5,3)(6.5,3)
\qbezier(4.5,3)(4.7,3.1)(4.9,3.2)
\qbezier(4.5,3)(4.7,2.9)(4.9,2.8)

\put(10,0){
\qbezier(9,1)(8,3)(7,5)
\qbezier(7,1)(7.4,1.8)(7.8,2.6)
\qbezier(9,5)(8.6,4.2)(8.2,3.4)
}

\linethickness{2.5pt}
\qbezier(13.00,2.00)(13.00,3.00)(13.00,4.00)
\linethickness{0.6pt}
\qbezier(12.00,1.00)(12.50,1.50)(13.00,2.00)
\qbezier(14.00,1.00)(13.50,1.50)(13.00,2.00)
\qbezier(12.00,5.00)(12.50,4.50)(13.00,4.00)
\qbezier(14.00,5.00)(13.50,4.50)(13.00,4.00)

\qbezier(14.00,5.00)(13.80,4.90)(13.60,4.80)
\qbezier(14.00,5.00)(13.90,4.80)(13.80,4.60)

\qbezier(12.00,5.00)(12.20,4.90)(12.40,4.80)
\qbezier(12.00,5.00)(12.10,4.80)(12.20,4.60)

\qbezier(12.50,1.50)(12.30,1.40)(12.10,1.30)\qbezier(12.50,1.50)(12.40,1.30)(12.30,1.10)

\qbezier(13.50,1.50)(13.70,1.40)(13.90,1.30)\qbezier(13.50,1.50)(13.60,1.30)(13.70,1.10)
\qbezier(14.5,3)(15.5,3)(16.5,3)
\qbezier(14.5,3)(14.7,3.1)(14.9,3.2)
\qbezier(14.5,3)(14.7,2.9)(14.9,2.8)

\qbezier(19.5,3)(20.5,3)(21.5,3)
\qbezier(21.5,3)(21.3,3.1)(21.1,3.2)
\qbezier(21.5,3)(21.3,2.9)(21.1,2.8)
\put(-10,0){
\qbezier(19,5)(18,3)(17,1)
\qbezier(17,5)(17.4,4.2)(17.8,3.4)
\qbezier(19,1)(18.6,1.8)(18.2,2.6)

\qbezier(19,5)(19,4.75)(19,4.5)
\qbezier(19,5)(18.8,4.8)(18.6,4.6)
\qbezier(17,5)(17,4.75)(17,4.5)
\qbezier(17,5)(17.2,4.8)(17.4,4.6)
}
\put(10,0){
\qbezier(9,5)(9,4.75)(9,4.5)
\qbezier(9,5)(8.8,4.8)(8.6,4.6)
\qbezier(7,5)(7,4.75)(7,4.5)
\qbezier(7,5)(7.2,4.8)(7.4,4.6)
}

\qbezier(22.00,1.00)(22.00,3.00)(22.00,5.00)
\qbezier(21.80,4.60)(21.90,4.80)(22.00,5.00)
\qbezier(22.20,4.60)(22.10,4.80)(22.00,5.00)

\qbezier(24.00,1.00)(24.00,3.00)(24.00,5.00)
\qbezier(23.80,4.60)(23.90,4.80)(24.00,5.00)
\qbezier(24.20,4.60)(24.10,4.80)(24.00,5.00)
\small{
\put(4.70,3.50){0-res.}
\put(9.6,3.50){1-res.}
\put(14.7,3.50){0-res.}
\put(19.6,3.50){1-res.}
}
\footnotesize{
\put(1.2,5.2){$x_1$}
\put(6.2,5.2){$x_1$}
\put(11.2,5.2){$x_1$}
\put(16.2,5.2){$x_1$}
\put(21.2,5.2){$x_1$}

\put(1.2,0.6){$x_4$}
\put(6.2,0.6){$x_4$}
\put(11.2,0.6){$x_4$}
\put(16.2,0.6){$x_4$}
\put(21.2,0.6){$x_4$}

\put(4.2,5.2){$x_2$}
\put(9.2,5.2){$x_2$}
\put(14.2,5.2){$x_2$}
\put(19.2,5.2){$x_2$}
\put(24.2,5.2){$x_2$}

\put(4.2,0.6){$x_3$}
\put(9.2,0.6){$x_3$}
\put(14.2,0.6){$x_3$}
\put(19.2,0.6){$x_3$}
\put(24.2,0.6){$x_3$}

}
}
\put(2.5,-1.5){$\Gamma_0$}
\put(22.5,-1.5){$\Gamma_0$}
\put(12.5,-1.5){$\Gamma_1$}

\end{picture}
\end{center}
\bigskip
\bigskip
and we call the resolutions obtained $\Gamma_0$ and $\Gamma_1$, as written in the picture. To the resolution $\Gamma_0$ we assign the following complex: 
$$\C(\Gamma_0):\quad\quad\quad 0\longrightarrow R\{-1,1\}\xrightarrow{x_2-x_3} R\longrightarrow 0$$
and to the resolution $\Gamma_1$ with the wide edge we assign the complex: 
$$\C(\Gamma_1):\quad\quad\quad 0\longrightarrow R\{-1,3\}\xrightarrow{(x_2-x_3)(x_4-x_2)} R\longrightarrow 0.$$

Assume that there are no free circles in the diagram $D$.
If we resolve all the crossings of $D$ we obtain a trivalent graph with wide edges. There are $2^n$ such resolutions $\Gamma$ of $D$ and to each one we assign the tensor product of $\C(\Gamma_0)$ and $\C(\Gamma_1)$, over all crossings $c(D)$, depending
on the type of resolution that appeared. In this way we have obtained a complex $\C(\Gamma)$ and to each resolution $\Gamma$ we will  assign  its cohomology $H(\Gamma)=H(\C(\Gamma))$.  \\
\indent Like in \cite{KR2}, we can obtain that $H(\Gamma)$ categorif\mbox{}ies the relations (\ref{rel1})--(\ref{rel4}). For example, the relation
(\ref{rel1}) becomes 
$$H(\Gamma\cup \textrm{ unknot})\cong
(H(\Gamma)\otimes\Q[x_i])\oplus 
(H(\Gamma)\otimes\Q[x_i]\{-1,1\}),$$
where $x_i$ is the label assigned to the circle (unknot). Note that in all def\mbox{}initions only the dif\mbox{}ferences $x_i-x_j$ appear. Thus we can work with the smaller ring of polynomials $R''=\Q[x_2-x_1,\ldots,x_{2n}-x_1]$ instead of $R'$ (like in \cite{KR2}). \\
\indent If we have free circles in the the diagram $D$, we introduce new variable $y$, with $\deg y=(0,2)$, extend the ring of polynomials to $R[y]$ and replace $R$ by $R[y]$ in the complexes $\C(\Gamma_i)$, $i=0,1$. Finally to every free circle we assign the complex: 
$$0\longrightarrow R'[y]\{-1,1\}\xrightarrow{y} R'[y]\longrightarrow 0,$$
and we tensor these complexes with $\C(\Gamma)$. In this way we obtain good value of the unknot (\ref{rel0}),  i.e. $H(U)\cong \Q$.\\
\indent We again organize the $2^n$ total resolutions $\Gamma$ of the diagram $D$ in the same cubic complex as in the standard categorif\mbox{}ications. To each vertex of the cube (i.e. to each total resolution $\Gamma$) we assign the graded vector space $H(\Gamma)$.

\indent We will introduce the dif\mbox{}ferentials between those cohomology groups as the maps induced by the (grading preserving) homomorphisms between the corresponding complexes $\C(\Gamma)$.
 Since these complexes are built as the tensor products of $\C(\Gamma_0)$ and $\C(\Gamma_1)$ it is enough to specify the homomorphisms between these two complexes. For a positive crossing $c$ we def\mbox{}ine  the following complex of complexes: 
\begin{equation}\C_c:\quad\quad 0\longrightarrow\C(\Gamma_0)\{0,2\}\stackrel{\chi_0}{\longrightarrow}\C(\Gamma_1)\longrightarrow 0,\label{dif1}\end{equation}
where $\C(\Gamma_1)$ is in cohomological degree 0, and the map $\chi_0$
is given by \\

\begin{equation*}
\begin{CD}
0@>>>R\{-1,3\}@>{x_2-x_3}>> R\{0,2\}@>>>0\\
 @.  @V{1}VV     @VV{{x_4-x_2}}V @.   \\
0@>>> R\{-1,3\}@>{(x_2-x_3)(x_4-x_2)}>>R@>>> 0.\end{CD}
\end{equation*}\\


\indent For a negative crossing $c$ we def\mbox{}ine  the following complex of complexes: 
\begin{equation}\C_c:\quad\quad 0\longrightarrow\C(\Gamma_1)\{0,-2\}\stackrel{\chi_1}{\longrightarrow}\C(\Gamma_0)\{0,-2\}\longrightarrow 0,\label{dif2}\end{equation}
where $\C(\Gamma_1)$ is in cohomological degree 0, and the map $\chi_1$
is given by \\
\begin{equation*}
\begin{CD}
0@>>>R\{-1,1\}@>{(x_2-x_3)(x_4-x_2)}>> R\{0,-2\}@>>>0\\
 @.  @V{x_4-x_2}VV     @VV{1}V @.   \\
0@>>> R\{-1,-1\}@>{x_2-x_3}>>R\{0,-2\}@>>> 0.\end{CD}
\end{equation*}\\


\indent Def\mbox{}ine $\C(D)$ as the tensor product of $\C_c$ over all crossings $c(D)$. It is a complex built out of Koszul complexes $\C(\Gamma)$, over all the total resolutions $\Gamma$ of the diagram $D$, and the dif\mbox{}ferential preserves the bigrading of each term $\C^j(D)$. Every $\C^j(D)$ decomposes as a direct sum of contractible two-term complexes and its cohomology $H(\C^j(D))$, which is denoted by $\C H^j(D)$. The dif\mbox{}ferential induces the grading preserving maps $\delta$ from $\C H^j(D)$ to $\C H^{j+1}(D)$ and we denote the complex obtained in this way by $\C H(D)$. The cohomology $H(D)=H(\C H(D),\delta)$ is triply-graded: 
$$H(D)=\bigoplus_{j,k,l}{H_{k,l}^j(D)}.$$
Here $j$ is the cohomology degree, and $k$ and $l$ come from the internal bigrading of the chain groups.\\
\indent In complete analogy with \cite{KR2} we have that $H(D)$ does not depend on the choice of the braid presentation $D$ of a link $L$, up to  an overall shift in the triple grading. Also, as we saw since $H(\Gamma)$ categorif\mbox{}ies the relations (\ref{rel0})--(\ref{rel4}) and since the dif\mbox{}ferentials are induced by the grading preserving maps (\ref{dif1}) and (\ref{dif2}) which obviously categorify the relations (\ref{1}) and (\ref{2}), we have that the 
 bigraded Euler characteristic of $\C H(D)$ is equal to $F(D)$. Finally, by introducing half-integral shifts as in \cite{wu} (in order to compensate the powers of $\alpha$ from (\ref{xxx})) by:
$$\H(D)=H(D)[\omega(D)/2]\{-\omega(D)/2,-\omega(D)/2\},$$
we obtain a triply-graded homology theory, which does not depend on the choice of the braid presentation $D$ of a link $L$ and whose bi-graded Euler characteristic is equal to the two-variable HOMFLYPT polynomial of a link $L$.

\section{New categorif\mbox{}ications of the chromatic and  dichromatic polynomials for graphs}\label{newgraf}

\subsection{Introduction}

\indent In this section we will def\mbox{}ine a complex of doubly-graded modules whose doubly-graded Euler characteristic is equal to the whole two-variable dicromatic polynomial.  The idea is partially inspired by the categorif\mbox{}ication of HOMFLYPT polynomial described in the previous section. \\
\indent Also, we
give a new categorif\mbox{}ication of the chromatic polynomial for graphs. We do this here dif\mbox{}ferently to \cite{gr}. We will def\mbox{}ine the chain groups (the modules corresponding to the vertices of the cube of resolutions) as the cohomologies of certain chain complexes.\\
\indent A graph $G$ is specif\mbox{}ied by a set of vertices $V(G)$ and a set of
edges $E(G)$. If $e$ is an arbitrary edge of the graph $G$,
then by $G-e$ we denote the graph $G$ with the edge $e$ deleted, and
by $G/e$ the graph obtained by contracting edge $e$ (i.e. by
identifying the vertices incident to $e$ and deleting $e$). 

\subsection{The chromatic polynomial}

 If $q$ is a positive integer, the chromatic polynomial $P_G(q)$ is def\mbox{}ined as the number of ways to color the vertices of $G$ by using at most $q$ colors, such that every two vertices which are connected by an edge receive a dif\mbox{}ferent color.
It is well-known that the chromatic polynomial can be def\mbox{}ined equivalently by the following two  axioms:
\begin{eqnarray*}
(C1)&\quad P_G=P_{G-e}-P_{G/e},\\
(C2)&\quad P_{N_k}=q^k,
\end{eqnarray*}
where $N_k$ is the graph with $k$ vertices and no edges. By using these axioms we can obviously extend the domain of the polynomial to the set of complex numbers, and, furthermore, instead of $q$ in the axiom $(C2)$   we will put $1/(1-q)$, with $|q|<1$.\\
\indent By repeated use of $(C1)$ (which is the famous deletion-contraction rule) we will obtain the value of the chromatic polynomial as a sum of  contributions from all spanning subgraphs of $G$ (subgraphs that contain all vertices of $G$), which we will call states.
Furthermore, if for each subset $s\subset E(G)$ we denote by $[G:s]$ the graph whose set of vertices is $V(G)$ and set of edges is $s$, then the contribution of the graph $[G:s]$ is $(-1)^{|s|}(1-q)^{-k(s)}$, where $|s|$ is the number of elements of $s$ and $k(s)$ is the number of connected components of $[G:s]$. Hence, we obtain the expression:
$$P_G(q)=\sum_{s\subset E(G)}{(-1)^{|s|}(1-q)^{-k(s)}}=\sum_{i\ge 0}{(-1)^i \sum_{s\subset E(G),|s|=i}{(1-q)^{-k(s)}}},$$
which is called the state-sum expansion of the polynomial $P_G(q)$.\\
 \indent In  subsection \ref{iiglavni} we will def\mbox{}ine a graded chain complex of modules $\C(G)$ whose graded Euler characteristic is equal to $P_G(q)$.\\

\subsection{The dichromatic polynomial}

The dichromatic polynomial $P_G(q,v)$ of the graph $G$ is a two-variable generalization of the chromatic polynomial given by the following two axioms:
\begin{eqnarray*}
(D1)&\quad P_G=P_{G-e}-qP_{G/e},\\
(D2)&\quad P_{N_k}=v^k,
\end{eqnarray*}
where $N_k$ is the graph with $k$ vertices and no edges.\\
\indent From $(D1)$ we have a recursive expression for the dichromatic polynomial in terms of the value of the polynomial on graphs with a smaller number of edges. Indeed, as in the case of the chromatic polynomial we obtain that the contribution of the state $[G:s]$ is $(-1)^{|s|}q^{|s|}v^{k(s)}$, where $|s|$ is the number of elements of $s$ and $k(s)$ is the number of connected components of $[G:s]$. Hence, we obtain the expression:
$$P_G(q,v)=\sum_{s\subset E(G)}{(-1)^{|s|}q^{|s|}v^{k(s)}}=\sum_{i\ge 0}{(-1)^i q^i\sum_{s\subset E(G),|s|=i}{v^{k(s)}}},$$
which is called the state-sum expansion of the polynomial $P_G(q,v)$. However, we will use a slightly dif\mbox{}ferent parametrization of the dichromatic polynomial, given by:
$$D_G(t,q)=(1+t^{-1}q)^{m}P_G(q,\frac{1+t^{-1}q}{1-q}),$$
where 
$m$ is the number of edges of the graph $G$.\\
\indent In subsection \ref{iiglavni2} we will def\mbox{}ine a chain complex $\D(G)$ of doubly graded modules whose doubly graded Euler characteristic is equal to $D_G(t,q)$.

\subsection{The categorif\mbox{}ication of the chromatic polynomial}\label{iiglavni}

Let $n$ denote the number of vertices of the graph $G$. Let $R$ be the ring of polynomials in $n$ variables over $\Q$, i.e. $R=\Q[x_1,\ldots,x_n]$. We introduce a grading in $R$, by giving the degree 1 to every $x_i$. Order the set of vertices of $G$ and to the $i$-th  vertex assign the variable $x_i$. Finally, to every edge $e\in E(G)$, whose endpoints are the vertices $i_e$ and $j_e$, assign the monomial $m_e=x_{i_e}-x_{j_e}$ (the ambiguity of the sign does not af\mbox{}fect the later construction).

\subsubsection{The cubic complex construction}
\indent Let $s\subset E(G)$ be a subset of the set of edges of $G$, and let $[G:s]$ be the corresponding state in the resolution of a graph $G$. Def\mbox{}ine the ideal $I_s$ as the ideal generated by the monomials $m_e$, for all edges $e\in s$. Finally, to the state $[G:s]$ assign the module $R_s=R/{I_s}$.\\
\begin{proposition}{\label{prop1}} The quantum graded dimension of $R_s$ is equal to
$(1-q)^{-k(s)}$, where $k(s)$ denotes the number of connected components of $[G:s]$.
\end{proposition}

\textbf{Proof:}\\
\indent Let $i$ and $j$ be two arbitrary vertices of $G$. They obviously belong to the same connected component of $[G:s]$ if and only if there exist a sequence of edges belonging to $s$ which connects $i$ and $j$, which obviously happens if and only if $x_i-x_j$ belongs to $I_s$. Hence, all the variables corresponding to the vertices from the same component, must be the same in $R_s$. In other words, $R_s$ is isomorphic to the ring of polynomials (over $\Q$) in $k(s)$ variables, and hence:
$$q\dim R_s={\left(\sum_{i\ge 0}{q^i}\right)}^{k(s)}={(1-q)}^{-k(s)}.$$ \kraj
\\
\\
\indent Denote by $m$ the number of edges of $G$, and f\mbox{}ix anordering on the set $E(G)$, denoted by $(e_1,\ldots,e_m)$. Now we will def\mbox{}ine the chain complex $\C$ in a standard way, by ''summing over columns" of our cubic complex:
for each $i$, with $0\le i \le m$, we will def\mbox{}ine the $i$-th chain group, ${\C^i(G)}$ as the direct sum of $R_s$, over all $s\subset E(G)$, such that $|s|=i$. \\
\indent Now, let us turn to the dif\mbox{}ferential. We will def\mbox{}ine the map $d^i$ from ${\C^i(G)}$ to ${\C^{i+1}(G)}$ as a sum of maps between the direct summands of the chain groups. The only nonzero maps are the maps from $R_s$ to $R_{s\cup \{e\}}$, with $e \notin s$ (which are exactly the ones that correspond to the edges of the cube), and we set them (up to a sign) to be the identity (i.e. the map that sends $f+I_s$ to $f+I_{s\cup \{e\}}$ for every $f\in R$). \\
\indent We now introduce signs in a standard way in order to make the cube anticommutative, and hence to make the square of the dif\mbox{}ferential equal to 0. Namely, we put minus signs exactly for those maps $R_s\to R_{s\cup \{e\}}$, with an odd number of edges in $s$ which are ordered before $e$.  \\
\indent In this way we have obtained a chain complex, $\C(G)$, of  graded $R$-modules with grading preserving dif\mbox{}ferential. Its homology groups obviously don't depend of the ordering of the vertices, and also don't depend of the ordering of the edges of $G$ (like in \cite{gr}, section 2.2.3), and hence we obtain
\begin{theorem}  
The homology groups of the chain complex $\C(G)$ are invariants of the graph $G$, and the graded Euler characteristic of $\C(G)$ is equal to the chromatic polynomial $P_G(q)$.
\end{theorem}

\subsubsection{Alternative description}

Now we will give an equivalent def\mbox{}inition of the chain complex $\C(G)$ in terms of Koszul complexes.\\
\indent To each edge $e\in E(G)$ we assign two complexes, $\C_{e-}$ and $\C_{e+}$ def\mbox{}ined in the following way:
$$\C_{e-}: \quad 0 \longrightarrow R \buildrel{0}\over\longrightarrow R \longrightarrow 0,$$
$$\C_{e+}: \quad 0 \longrightarrow R \buildrel{x_i-x_j}\over\longrightarrow R \longrightarrow 0,$$
where $i$ and $j$ are the vertices of $G$ which are the endpoints of the edge $e$.
Now, to every subset $s\subset E(G)$ we assign a complex $\C_s$ which is the tensor product of $\C_{e\pm}$, where we take $+$ if $e\in s$ and $-$ if $e \notin s$. Finally, to the state $[G:s]$ we assign the cohomology of $\C_s$ at the rightmost position.\\
\indent To build the dif\mbox{}ferentials, we introduce the (grading preserving) maps $d_e$, as the maps induced on the cohomology by the following homomorphism from $\C_{e-}$ to $\C_{e+}$:\\
\begin{equation}
\begin{CD}
0@>>>R@>{0}>> R@>>>0\\
 @.  @V{0}VV     @VV{1}V @.  \\
0@>>> R@>{x_i-x_j}>>R@>>> 0.\end{CD}
\label{dif}\end{equation}\\

Here we put the upper row in cohomological degree 0, and the lower one in cohomological degree 1.\\
\indent Now, in order to def\mbox{}ine the dif\mbox{}ferentials, just tensor all the
chain complexes and maps between them from (\ref{dif}) over all edges $e$ of $E(G)$. If we take the cohomology only at the rightmost position in  each ``horizontal" complex (the ones in the same cohomological degree with respect to the def\mbox{}inition after (\ref{dif})), and as the dif\mbox{}ferentials are the induced maps between them, we obtain a complex  $\C'(G)$ which is isomorphic to the complex $\C(G)$ from the previous subsection.

\subsection{The categorif\mbox{}ication of the dichromatic polynomial}\label{iiglavni2}

In order to categorify the dichromatic polynomial we will have to introduce a new grading direction, and we will use the whole Koszul complex (actually a slightly modif\mbox{}ied one) that we have used in the previous subsection.\\
\indent We order the vertices of $G$, and to the $i$-th one ($1\le i \le n=\sharp V(G))$, we assign the variable $x_i$. We def\mbox{}ine the bidegree of all $x_i$ as (0,1). Def\mbox{}ine the bigraded ring $R$ by $R=\Q[x_1,\ldots,x_n]$, where we put the f\mbox{}ield $\Q$ in bidegree (0,0).\\
\indent To each edge $e$, such that its endpoints are the $i$-th and $j$-th vertex, we can associate two resolutions of the graph $G$: the f\mbox{}irst one with the edge $e$ contracted (i.e. when we identify the vertices $i$ and $j$), and the second one with the edge $e$ deleted. To the f\mbox{}irst resolution we assign the following complex (denoted by $\D(e+)$):
$$\D(e+):\qquad 0 \longrightarrow R\{-1,1\} \stackrel{x_i-x_j}{\longrightarrow} R \longrightarrow 0,$$
and to the second one we assign the complex $\D(e-)$ given by:
$$\D(e-):\qquad 0 \longrightarrow R\{-1,1\} \stackrel{0}{\longrightarrow} R \longrightarrow 0.$$
Let $s$ be an arbitrary subset of $E(G)$ and let $[G:s]$ be the corresponding state of $G$. Then, to that state we assign the (Koszul) complex $\D'(s)$ of bigraded $R$-modules obtained by tensoring the complexes $\D(e+)$, where $e$ runs over all edges in $s$, and $\D(f-)$, where $f$ runs over all edges in $E(G)\setminus s$. We denote its (bigraded) cohomology by $H'(s)$ (the direct sum of the cohomology groups of $\D'(s)$).
\begin{proposition}
The quantum bigraded dimension of $H'(s)$ is equal to:
$${(1+t^{-1}q)}^{m-n}{\left(\frac{1+t^{-1}q}{1-q}\right)}^{k(s)},$$
where $k(s)$ is the number of connected components of $[G:s]$, and $n$ and $m$ are the number of vertices and edges of $G$, respectively. 
\end{proposition}
\textbf{Proof:}\\
\indent Like in the proof of Proposition \ref{prop1} we obtain that the cohomology at the rightmost position of $\D'(s)$ is isomorphic to the ring of polynomials in $k(s)$ variables. However, here we will also have the cohomology at the leftmost position in each of the $\D(e\pm)$, which is isomorphic to the same ring of polynomials in $k(s)$ variables, but shifted by the bidegree $\{-1,1\}$, for all $\D(e-)$ and for a certain number of the $D(e+)$. We will show by induction on $|s|$ that the total number of such $e$'s, denoted by $c(s)$, is equal to $k(s)-n+m$.\\
\indent If $|s|=0$ then we have the tensor product of $m$ complexes with all the mappings equal to zero, and hence we have that the number of edges which contribute with nontrivial cohomology at the leftmost position is equal to $m=k(s)-n+m$ (note that in this case $k(s)=n$). Now suppose that the formula is true for some subset $s$ and consider the state $[G:{(s\cup e)}]$ with $e \in E(G) \setminus s$. Denote the endpoints of $e$ by $i$ and $j$, and denote $s'=s\cup e$. This means that  $\D'(s')$ is formed by the tensor product of the same complexes as $\D'(s)$ with $\D(e+)$ instead of $\D(e-)$. Now, $\D(e+)$ will have nontrivial cohomology at the leftmost position if and only if $x_i-x_j$ belongs to the ideal generated by the monomials def\mbox{}ined by the edges of $s$, i.e. if and only if the vertices $i$ and $j$ belong to the same connected component of $[G:s]$. In other words, we have $c(s')=c(s)$ if $k(s')=k(s)$ and $c(s')=c(s)-1$ if $k(s')=k(s)-1$. So we have $c(s)=k(s)-n+m$ as we wanted to prove.\\
\indent Hence the total bigraded dimension of $H'(s)$ is equal to:
$${(1+t^{-1}q)}^{k(s)-n+m}{(1-q)}^{-k(s)}.$$ \kraj

Furthermore, to every vertex $v$ of the graph $G$, we assign the same complex as $\D(e-)$. Now, if we tensor these complexes over all vertices of $G$ and tensor the complex obtained with $\D'(s)$, we obtain the complex $\D(s)$. We denote the cohomology of $\D(s)$ by $H(s)$, and that is the space that we will assign to the state $[G:s]$. Obviously, we have:
$$q\dim H(s)={(1+t^{-1}q)}^m{\left(\frac{1+t^{-1}q}{1-q}\right)}^{k(s)}. $$

\indent In order to introduce the dif\mbox{}ferentials between the cohomologies $H(s)$, we will def\mbox{}ine the (grading preserving) homomorphism $d(e)$ from $\D(e+)\{0,1\}$ to $\D(e-)$, and then for the dif\mbox{}ferentials we take the induced mappings on the cohomology. We def\mbox{}ine $d(e)$ by:\\
\begin{equation*}
\begin{CD}
0@>>>R\{-1,2\}@>{x_i-x_j}>> R\{0,1\}@>>>0\\
 @.  @V{x_i-x_j}VV     @VV{0}V @.   \\
0@>>> R\{-1,1\}@>{0}>>R@>>> 0.\end{CD}
\end{equation*}\\

We put the upper row in cohomological degree -1, and the lower row in cohomological degree 0. We will denote this complex of complexes 
by $\D_e$. \\
\indent For the graph $G$ def\mbox{}ine the complex of (Koszul) complexes by tensoring $\D_e$ over all edges $e$ of $G$. By taking the cohomology $H^j(G)$, $-m\le j\le 0$,  in each ``horizontal" complex and def\mbox{}ining the dif\mbox{}ferentials between them to be the ones induced by the tensor product of $d(e)$'s, we obtain a triply graded  complex $\D(G)$.\\
\indent From the def\mbox{}inition we have

\begin{theorem}
The homotopy class of the complex $\D(G)$ is an invariant of the graph $G$ whose bigraded Euler characteristic is equal to the  dichromatic polynomial $D_G(t,q)$ of the graph $G$.
\end{theorem}

\end{document}